\documentclass[reqno]{article}%
\usepackage{amsmath}
\usepackage{graphicx}
\usepackage{amsfonts}
\usepackage{amssymb}%
\newtheorem{theorem}{Theorem}

\newtheorem{corollary}[theorem]{Corollary}

\newtheorem{example}[theorem]{Example}

\newenvironment{proof}[1][Proof]{\noindent{\textbf {#1}  }}  {\hfill$\Box$\bigskip}

\begin{document}

\title{Books in graphs}
\author{B\'{e}la Bollob\'{a}s\thanks{Department of Mathematical Sciences, University
of Memphis, Memphis TN 38152, USA} \thanks{Trinity College, Cambridge CB2 1TQ,
UK} \thanks{Research supported in part by NSF grant DSM 9971788 and DARPA
grant F33615-01-C-1900.} \ and Vladimir Nikiforov$^{\dag}$}
\maketitle

\begin{abstract}
A set of $q$ triangles sharing a common edge is called a book of size $q.$ We
write $\beta\left(  n,m\right)  $ for the the maximal $q$ such that every
graph $G\left(  n,m\right)  $ contains a book of size $q$. In this note

1) we compute $\beta\left(  n,cn^{2}\right)  $ for infinitely many values of
$c$ with $1/4<c<1/3$,

2) we show that if $m\geq\left(  1/4-\alpha\right)  n^{2}$ with $0<\alpha
<17^{-3},$ and $G$ has no book of size at least $\left(  1/6-2\alpha
^{1/3}\right)  n$ then $G$ contains an induced bipartite graph $G_{1}$ of
order at least $\left(  1-\alpha^{1/3}\right)  n$ and minimal degree
\[
\delta\left(  G_{1}\right)  \geq\left(  \frac{1}{2}-4\alpha^{1/3}\right)  n,
\]

3) we apply the latter result to answer two questions of Erd\H{o}s concerning
the booksize of graphs $G\left(  n,n^{2}/4-f\left(  n\right)  n\right)  $
every edge of which is contained in a triangle, and $0<f\left(  n\right)
<n^{2/5-\varepsilon}.$

\end{abstract}

\section{Introduction}

Our notation and terminology are standard (see, e.g., \cite{Bol}). Thus,
$G\left(  n,m\right)  $ is a graph of order $n$ and size $m;$ for a graph $G$
and a vertex $u\in V\left(  G\right)  $ we write $\Gamma\left(  u\right)  $
for the set of vertices adjacent to $u;$ $d_{G}\left(  u\right)  =\left\vert
\Gamma\left(  u\right)  \right\vert $ is the degree of $u;$ we write $d\left(
u\right)  $ instead of $d_{G}\left(  u\right)  $ if the graph $G$ is implicit.
However, somewhat unusually, we set $\widehat{d}\left(  U\right)  =\left\vert
\cap_{x\in U}\Gamma\left(  x\right)  \right\vert $. Unless explicitly stated,
all graphs are assumed to be defined on the vertex set $\left[  n\right]
=\left\{  1,2,...n\right\}  .$ Also, $k_{s}\left(  G\right)  $ is the number
of $s$-cliques of $G.$

In 1962 Erd\H{o}s \cite{Erd1} initiated the study of books in graphs. A
\emph{book }of size $q$ consists of $q$ triangles sharing a common edge. We
write $bk\left(  G\right)  $ for the size of the largest book in a graph $G$
and call it the \emph{booksize }of\emph{ }$G.$ Since 1962 books have attracted
considerable attention both in extremal graph theory (see, e.g., \cite{KhNi},
\cite{EFR}, and \cite{EFG}) and in Ramsey graph theory (see, e.g., \cite{RS},
\cite{FRS}, and \cite{NR}).

Erd\H{o}s, Faudree and Rousseau defined in \cite{EFR} the function
\[
\beta\left(  n,m\right)  =\min\left\{  bk\left(  G\right)  |\text{ }G=G\left(
n,m\right)  \right\}  .
\]
Our aim in this paper the study of the function $\beta\left(  n,m\right)  $
and its variants. We shall prove a technical inequality about booksizes that
we shall use to give bounds on $\beta\left(  n,m\right)  $ and answer two
questions of Erd\H{o}s.

The paper is organized as follows: in section 2 we use a counting argument of
Khad\v{z}iivanov and Nikiforov \cite{KhNi} to prove a bound on $\beta\left(
n,m\right)  $ in terms of the degree sequence and other graph parameters. In
particular, this result implies that $\beta\left(  n,\left\lfloor
n^{2}/4\right\rfloor +1\right)  >n/6,$ as conjectured by Erd\H{o}s and proved
by Edwards \cite{Edw}. In addition, we determine $\beta\left(  n,cn^{2}%
\right)  $ for infinitely many values of $c$ with $1/4<c<1/3$. In section 3 we
prove that a graph $G\left(  n,\left(  1/4-\alpha\right)  n^{2}\right)  $ with
$0<\alpha<17^{-3}$ either has a book of size about $n/6$ or has a large
induced bipartite graph with minimal degree close to $n/2.$ In the last
section we make use of this structural property to answer two questions of
Erd\H{o}s concerning the booksize of graphs $G\left(  n,n^{2}/4-f\left(
n\right)  n\right)  ,$ every edge of which is contained in a triangle and
$0<f\left(  n\right)  \leq n^{2/5-\varepsilon}$.

\section{A lower bound on the booksize of a graph}

In 1962 Erd\H{o}s \cite{Erd1} conjectured that the booksize of a graph $G$ of
order $n$ and size greater than $\left\lfloor n^{2}/4\right\rfloor $ is at
least $\left\lfloor n/6\right\rfloor ,$ i.e., $\beta\left(  n,\left\lfloor
n^{2}/4\right\rfloor +1\right)  \geq n/6$. This was proved by Edwards in an
unpublished manuscript \cite{Edw} and independently by Khad\v{z}iivanov and
Nikiforov in \cite{KhNi}.

For $r\geq3$ and $0\leq j<r,$ we write $K_{r}^{(j)}$ for the graph consisting
of a complete graph $K_{r-1}$ and an additional vertex joined to precisely
$r-j-1$ vertices of the $K_{r-1}.$ We denote by $k_{r}^{(j)}\left(  G\right)
$ the number of induced subgraphs of $G$ that are isomorphic to $K_{r}^{(j)},$
e. g., $k_{4}^{\left(  3\right)  }\left(  G\right)  $ is the number of induced
subgraphs of $G$ that are isomorphic to a triangle with an isolated vertex.

\begin{theorem}
\label{thBT} Let $G=G\left(  n,m\right)  $ be a graph with degree sequence
$d\left(  1\right)  ,...,d\left(  n\right)  .$ Then,
\[
\left(  6k_{3}\left(  G\right)  -\sum_{i=1}^{n}d^{2}\left(  i\right)
+nm\right)  bk\left(  G\right)  \geq nk_{3}\left(  G\right)  +8k_{4}\left(
G\right)  +2k_{4}^{\left(  3\right)  }\left(  G\right)  .
\]

\end{theorem}

\begin{proof}
In the proof we use some arguments from \cite{KhNi}. Set $\beta=bk\left(
G\right)  .$ Clearly $G$ contains exactly $\left(  n-3\right)  k_{3}\left(
G\right)  $ pairs $\left(  v,T\right)  $ where $v\in V\left(  G\right)  $
\ and $T$ is a triangle in $G.$ Also, a $K_{4}$ subgraph of $G$ contains
exactly $4$ such pairs; a $K_{4}^{\left(  j\right)  }$ subgraph contains two
such pairs for $j=1,$ and one such pair for $j=2$ and $3.$ Therefore,
\begin{equation}
\left(  n-3\right)  k_{3}\left(  G\right)  =4k_{4}\left(  G\right)
+2k_{4}^{\left(  1\right)  }\left(  G\right)  +k_{4}^{\left(  2\right)
}\left(  G\right)  +k_{4}^{\left(  3\right)  }\left(  G\right)  . \label{eq1}%
\end{equation}
We have%
\[
\sum_{\left(  i,j\right)  \in E\left(  G\right)  }\binom{\widehat{d}\left(
ij\right)  }{2}=6k_{4}\left(  G\right)  +k_{4}^{\left(  1\right)  }\left(
G\right)  ,
\]
yielding
\[
\sum_{\left(  i,j\right)  \in E\left(  G\right)  }\left(  \widehat{d}%
^{2}\left(  ij\right)  -\widehat{d}\left(  ij\right)  \right)  =12k_{4}\left(
G\right)  +2k_{4}^{\left(  1\right)  }\left(  G\right)  .
\]
Since
\begin{equation}
\sum_{\left(  i,j\right)  \in E\left(  G\right)  }\widehat{d}\left(
ij\right)  =3k_{3}\left(  G\right)  , \label{eq11}%
\end{equation}
we see that
\[
\sum_{\left(  i,j\right)  \in E\left(  G\right)  }\widehat{d}^{2}\left(
ij\right)  =12k_{4}\left(  G\right)  +2k_{4}^{\left(  1\right)  }\left(
G\right)  +3k_{3}\left(  G\right)  .
\]
Subtracting (\ref{eq1}) from the last equality and rearranging the terms, we
obtain
\begin{equation}
nk_{3}\left(  G\right)  =\sum_{\left(  i,j\right)  \in E\left(  G\right)
}\widehat{d}^{2}\left(  ij\right)  -8k_{4}\left(  G\right)  +k_{4}^{\left(
2\right)  }\left(  G\right)  +k_{4}^{\left(  3\right)  }\left(  G\right)  .
\label{eq2}%
\end{equation}
Next we shall eliminate the term $k_{4}^{\left(  2\right)  }\left(  G\right)
$ from (\ref{eq2}). For every $i\in V\left(  G\right)  $ set $\Gamma^{\prime
}\left(  i\right)  =V\left(  G\right)  \backslash\Gamma\left(  i\right)  .$
The sum $\sum_{ij\in E\left(  G\right)  }\widehat{d}\left(  ij\right)
\left\vert \Gamma^{\prime}\left(  i\right)  \cap\Gamma^{\prime}\left(
j\right)  \right\vert $ counts each $K_{4}^{(2)}$ once and each $K_{4}^{(3)}$
three times, so
\begin{equation}
\sum_{\left(  i,j\right)  \in E\left(  G\right)  }\widehat{d}\left(
ij\right)  \left\vert \Gamma^{\prime}\left(  i\right)  \cap\Gamma^{\prime
}\left(  j\right)  \right\vert =k_{4}^{\left(  2\right)  }\left(  G\right)
+3k_{4}^{\left(  3\right)  }\left(  G\right)  . \label{four}%
\end{equation}
Subtracting (\ref{four}) from (\ref{eq2}), we see that
\begin{align}
nk_{3}\left(  G\right)   &  =\sum_{\left(  i,j\right)  \in E\left(  G\right)
}\widehat{d}^{2}\left(  ij\right)  +\sum_{\left(  i,j\right)  \in E\left(
G\right)  }\widehat{d}\left(  ij\right)  \left\vert \Gamma^{\prime}\left(
i\right)  \cap\Gamma^{\prime}\left(  j\right)  \right\vert -8k_{4}\left(
G\right)  -2k_{4}^{\left(  3\right)  }\left(  G\right) \nonumber\\
&  =\sum_{\left(  i,j\right)  \in E\left(  G\right)  }\widehat{d}\left(
ij\right)  \left(  \widehat{d}\left(  ij\right)  +\left\vert \Gamma^{\prime
}\left(  i\right)  \cap\Gamma^{\prime}\left(  j\right)  \right\vert \right)
-8k_{4}\left(  G\right)  -2k_{4}^{\left(  3\right)  }\left(  G\right)  .
\label{five}%
\end{align}
Noting that $\widehat{d}\left(  ij\right)  \leq\beta$ for every edge $\left(
i,j\right)  $ and recalling (\ref{eq11}), inequality (\ref{five}) implies
that
\begin{align}
nk_{3}\left(  G\right)   &  \leq\beta\sum_{\left(  i,j\right)  \in E\left(
G\right)  }\left(  \widehat{d}\left(  ij\right)  +\left\vert \Gamma^{\prime
}\left(  i\right)  \cap\Gamma^{\prime}\left(  j\right)  \right\vert \right)
-8k_{4}\left(  G\right)  -2k_{4}^{\left(  3\right)  }\left(  G\right)
\nonumber\\
&  =\beta\left(  3k_{3}\left(  G\right)  +\sum_{\left(  i,j\right)  \in
E\left(  G\right)  }\left\vert \Gamma^{\prime}\left(  i\right)  \cap
\Gamma^{\prime}\left(  j\right)  \right\vert \right)  -8k_{4}\left(  G\right)
-2k_{4}^{\left(  3\right)  }\left(  G\right)  \label{eq3}%
\end{align}
Since
\[
\left\vert \Gamma^{\prime}\left(  i\right)  \cap\Gamma^{\prime}\left(
j\right)  \right\vert =n-d\left(  i\right)  -d\left(  j\right)  +\widehat
{d}\left(  ij\right)  .
\]
we find that
\begin{align*}
\sum_{\left(  i,j\right)  \in E\left(  G\right)  }\left\vert \Gamma^{\prime
}\left(  i\right)  \cap\Gamma^{\prime}\left(  j\right)  \right\vert  &
=\sum_{\left(  i,j\right)  \in E\left(  G\right)  }\left(  n-d\left(
i\right)  -d\left(  j\right)  +\widehat{d}\left(  ij\right)  \right) \\
&  =3k_{3}\left(  G\right)  +nm-\sum_{i=1}^{n}d^{2}\left(  i\right)  .
\end{align*}
Putting this into (\ref{eq3}) we see that
\[
nk_{3}\left(  G\right)  +8k_{4}\left(  G\right)  +2k_{4}^{\left(  3\right)
}\left(  G\right)  \leq6\beta k_{3}\left(  G\right)  +\beta\left(  -\sum
_{i=1}^{n}d^{2}\left(  i\right)  +nm\right)  ,
\]
as claimed.
\end{proof}

The following corollary is due to Edwards \cite{Edw}.

\begin{corollary}
\label{cobs} For every graph $G=$ $G\left(  n,m\right)  $ with $m>n^{2}/4$%
\begin{equation}
bk\left(  G\right)  \geq\frac{2m}{n}-\frac{n}{3}. \label{lobs}%
\end{equation}

\end{corollary}

\begin{proof}
With $\beta=bk\left(  G\right)  $, Theorem \ref{thBT} implies that
\[
\left(  6k_{3}\left(  G\right)  -\sum_{i=1}^{n}d^{2}\left(  i\right)
+nm\right)  \beta\geq nk_{3}\left(  G\right)  +8k_{4}\left(  G\right)
+2k_{4}^{\left(  3\right)  }\left(  G\right)  \geq nk_{3}\left(  G\right)  ,
\]
and so
\begin{equation}
\left(  6\beta-n\right)  k_{3}\left(  G\right)  \geq\beta\left(  \sum
_{i=1}^{n}d^{2}\left(  i\right)  -nm\right)  . \label{six}%
\end{equation}
Since $\sum_{i=1}^{n}d\left(  i\right)  =2m$, we have
\begin{equation}
\sum_{i=1}^{n}d^{2}\left(  i\right)  \geq\frac{4m^{2}}{n}>nm; \label{seven}%
\end{equation}
in particular,
\[
\sum_{i=1}^{n}d^{2}\left(  i\right)  -nm>0.
\]
Hence, (\ref{six}) implies that $6\beta>n$. Furthermore, as $3k_{3}\left(
G\right)  \leq\beta m$, we see from (\ref{six}) and (\ref{seven}) that
\[
\frac{1}{3}\left(  6\beta-n\right)  \beta m\geq\beta\left(  \frac{4m^{2}}%
{n}-nm\right)  ,
\]
implying (\ref{lobs}).
\end{proof}

As a consequence of Corollary \ref{cobs} we easily obtain the following bound.

\begin{corollary}
\label{co1} For every graph $G\left(  n,\left\lfloor n^{2}/4\right\rfloor
+1\right)  $ we have $bk\left(  G\right)  >n/6$.\hfill$\square$
\end{corollary}

The graph $H_{s,t}$ below, constructed by Erd\H{o}s, Faudree and Rousseau in
\cite{EFR}, shows that the bound in Corollary \ref{cobs} is essentially best possible.

\begin{example}
\label{examp1} Let $t\geq1,$ $s>3$ be fixed integers. Partition the vertex set
$V=\left[  n\right]  $ with $n=3st$ into $3s$ sets $V_{ij}$ $\left(
i\in\left[  3\right]  ,\text{ }j\in\left[  r\right]  \right)  $ of cardinality
$t.$ Join two vertices $v\in V_{ij}$ and $u\in V_{kl}$ iff $i\neq k$ and
$j\neq l.$
\end{example}

By straightforward counting we see that
\[
e\left(  H_{s,t}\right)  =3s\left(  s-1\right)  t^{2}=3s\left(  s-1\right)
\left(  \frac{n}{3s}\right)  ^{2}=\frac{s-1}{3s}n^{2},
\]
and
\[
bk\left(  H_{s,t}\right)  =\left(  s-2\right)  t=\frac{\left(  s-2\right)
n}{3s}.
\]
On the other hand, from Corollary \ref{cobs},\ we have%
\[
bk\left(  H_{s,t}\right)  \geq\frac{2e\left(  H_{s,t}\right)  }{n}-\frac{n}%
{3}=\frac{2\left(  s-1\right)  n}{3s}-\frac{n}{3}=\frac{\left(  s-2\right)
n}{3s},
\]
thus, the bound in Corollary \ref{cobs} is tight for $n,m$ with $3s|n$, $s>8$,
and $m=\left(  s-1\right)  n^{2}/3s.$

A different extremal graph (\cite{Edw}, \cite{KhNi}) is defined as follows.

\begin{example}
\label{examp2} Select $6$ disjoint sets $A_{11},A_{12},A_{13},A_{21}%
,A_{22},A_{23}$ with $\left|  A_{11}\right|  =\left|  A_{12}\right|  =\left|
A_{13}\right|  =k-1$ and $\left|  A_{21}\right|  =\left|  A_{22}\right|
=\left|  A_{23}\right|  =k+1.$ Set $V\left(  G\right)  $ to be the union of
all these sets. For every $1\leq j<k\leq3$ join every vertex of $A_{ij}$ to
every vertex of $A_{ik}$ and for $j=1,2,3$ join every vertex of $A_{1j}$ to
every vertex of $A_{2j}.$
\end{example}

It is easy to check that the resulting graph has $n=6k$ vertices,
$9k^{2}+3>n^{2}/4$ edges and its booksize is precisely $k+1=n/6+1.$

\section{A stability theorem for graphs without large books}

In this section we give a structural property of graphs having substantial
size and whose booksize is small.

In \cite{AES} Andr\'{a}sfai, Erd\H{o}s and S\'{o}s proved that if $G$ is a
$K_{r+1}$-free graph of order $n$ with minimal degree
\[
\delta\left(  G\right)  >\left(  1-\frac{3}{3r-1}\right)  n
\]
then $G$ is $r$-chromatic. We shall use this theorem to obtain a structural
result related to the stability theorems of Simonovits (see, e. g., \cite{Si1}).

\begin{theorem}
\label{stab} \label{th2} For every $\alpha$ with $0<\alpha<10^{-5}$ and every
graph $G=G\left(  n,m\right)  $ with
\begin{equation}
m\geq\left(  \frac{1}{4}-\alpha\right)  n^{2} \label{cond1}%
\end{equation}
either
\begin{equation}
bk\left(  G\right)  >\left(  \frac{1}{6}-2\alpha^{1/3}\right)  n \label{prop1}%
\end{equation}
or $G$ contains an induced bipartite graph $G_{1}$ of order at least $\left(
1-\alpha^{1/3}\right)  n$ and with minimal degree
\begin{equation}
\label{prop2}\delta\left(  G_{1}\right)  \geq\left(  \frac{1}{2}-4\alpha
^{1/3}\right)  n.
\end{equation}

\end{theorem}

\begin{proof}
If $m>n^{2}/4$ then Corollary \ref{co1} implies that $bk\left(  G\right)
>n/6$, which is stronger than (\ref{prop1}), so we may assume that $m\leq
n^{2}/4.$ Furthermore, if $\sum_{i=1}^{n}d^{2}\left(  i\right)  >nm$ then
Theorem \ref{thBT} implies that
\[
\left(  6bk\left(  G\right)  -n\right)  k_{3}\left(  G\right)  >0,
\]
and so again $bk\left(  G\right)  >n/6$. Therefore, we may assume
\[
\sum_{i=1}^{n}d^{2}\left(  i\right)  \leq nm.
\]
Clearly, from (\ref{cond1}),
\[
\frac{4m^{2}}{n}\geq m\left(  n-4\alpha n\right)  =nm-4\alpha nm,
\]
and so,
\begin{equation}
\sum_{i=1}^{n}\left(  d\left(  i\right)  -\frac{2m}{n}\right)  ^{2}=\sum
_{i=1}^{n}d^{2}\left(  i\right)  -\frac{4m}{n}^{2}\leq4\alpha nm\leq\alpha
n^{3}. \label{in1}%
\end{equation}
Set $\varepsilon=\alpha^{1/3}$, $M=\{u\in V(G):\ d\left(  u\right)  <\frac
{2m}{n}-\varepsilon n\}$ and $G_{1}=G\left[  V\backslash M\right]  $. We claim
that $G_{1}$ has the required properties. First we show that its minimal
degree satisfies (\ref{prop2}). From (\ref{in1}),
\[
\left\vert M\right\vert \varepsilon^{2}n^{2}\leq\sum_{v\in M}\left(  d\left(
v\right)  -\frac{2m}{n}\right)  ^{2}<\sum_{i=1}^{n}\left(  d\left(  i\right)
-\frac{2m}{n}\right)  ^{2}\leq\alpha n^{3}.
\]
Hence, $\left\vert M\right\vert <\left(  \alpha/\varepsilon^{2}\right)
n=\alpha^{1/3}n$, i.e., $v\left(  G_{1}\right)  >\left(  1-\alpha
^{1/3}\right)  n.$ Also, for $v\in V\backslash M$, we have
\begin{align}
d_{G_{1}}\left(  v\right)   &  \geq d\left(  v\right)  -\left\vert
M\right\vert >\left(  \frac{2m}{n}-\varepsilon n\right)  -|M|=\frac{n}%
{2}-2\alpha n-\alpha^{1/3}n-|M|\nonumber\\
&  >\left(  \frac{1}{2}-2\alpha n-2\alpha^{1/3}\right)  n\geq\left(  \frac
{1}{2}-4\alpha^{1/3}\right)  n. \label{fourteen}%
\end{align}
All that remains to prove is that $G_{1}$ is bipartite. Suppose first that
$G_{1}$ contains a triangle with vertices $u,v,w$, say. Since
\[
n\geq d\left(  u\right)  +d\left(  v\right)  +d\left(  w\right)  -\widehat
{d}\left(  uv\right)  -\widehat{d}\left(  uw\right)  -\widehat{d}\left(
vw\right)
\]
we find that
\begin{align*}
\widehat{d}\left(  uv\right)  +\widehat{d}\left(  uw\right)  +\widehat
{d}\left(  vw\right)   &  \geq d\left(  u\right)  +d\left(  v\right)
+d\left(  w\right)  -n\\
&  \geq3\left(  \frac{1}{2}-\alpha-\sqrt[3]{\alpha}\right)  n-n.
\end{align*}
Thus,
\[
bk\left(  G\right)  \geq\left(  \frac{1}{6}-\alpha n-\alpha^{1/3}\right)
n\geq\left(  \frac{1}{6}-2\alpha^{1/3}\right)  n,
\]
and so (\ref{prop2}) holds. Finally, assume that $G_{1}$ is triangle-free.
Since $\alpha<10^{-5}$,
\[
\delta(G_{1})\geq\left(  \frac{1}{2}-4\alpha^{1/3}\right)  n>\frac{2}%
{5}v(G_{1}).
\]
Hence, the case $r=2$ of the theorem of Andr\'{a}sfai, Erd\H{o}s and S\'{o}s
mentioned above implies that $G_{1}$ is indeed bipartite, completing the proof
of Theorem \ref{stab}
\end{proof}

It is easily seen that if we are a little more careful in our proof of
$\delta(G_{1}) > v(G_{1})$ then the condition on $\alpha$ can be relaxed to
$0<\alpha<17^{-3}$.

\section{Two problems of Erd\H{o}s}

Erd\H{o}s and Rothschild suggested the study of the booksize of graphs in
which every edge is contained in a triangle. In \cite{Erd2} and \cite{Erd3}
Erd\H{o}s himself gave some results on such graphs. Suppose $f\left(
n\right)  $ is a fixed positive function of $n,$ and let $TG\left(
n,f\right)  $ be the set of all graphs $G=G\left(  n,m\right)  $ such that
every edge of $G$ is contained in a triangle and $m>\max\left\{
n^{2}/4-f\left(  n\right)  n,0\right\}  $. Set%
\[
\gamma\left(  n,f\right)  =\min\left\{  bk\left(  G\right)  \text{ }|\text{
}G\in TG\left(  n,f\right)  \right\}  .
\]
In \cite{Erd2}, p. 91, Erd\H{o}s proved that for every $c>0$ there exists some
$c_{1}>0$ such that
\[
\gamma\left(  n,c\right)  \geq c_{1}n
\]
for $n$ sufficiently large. Hence, setting%
\[
\underline{\lim}_{n\rightarrow\infty}\frac{\gamma\left(  n,c\right)  }%
{n}=\sigma\left(  c\right)  ,
\]
we see that for every $c>0,$ $\sigma\left(  c\right)  >0.$ Erd\H{o}s asked how
large $\sigma\left(  c\right)  $ is. Our next theorem gives an answer that is
asymptotically tight when $c$ tends to $0$.

\begin{theorem}
\label{th1} For every function $f\left(  n\right)  $ with $0<f\left(
n\right)  <n/4,$
\[
\gamma\left(  n,f\right)  >\frac{n}{12f\left(  n\right)  +6}.
\]

\end{theorem}

\begin{proof}
From Theorem \ref{th1} we have for $\beta=bk\left(  G\right)  $%
\[
\left(  6k_{3}\left(  G\right)  -\sum_{i=1}^{n}d^{2}\left(  i\right)
+nm\right)  \beta\geq nk_{3}\left(  G\right)  ,
\]
and hence,%
\[
\left(  6\beta-n\right)  k_{3}\left(  G\right)  \geq\beta\left(  \sum
_{i=1}^{n}d^{2}\left(  i\right)  -nm\right)  .
\]
From $\sum_{i=1}^{n}d\left(  i\right)  =2m$ we have $\sum_{i=1}^{n}%
d^{2}\left(  i\right)  \geq4m^{2}/n$ and thus,%
\[
\left(  6\beta-n\right)  k_{3}\left(  G\right)  \geq\beta\left(  \frac{4m^{2}%
}{n}-nm\right)  >-4f\left(  n\right)  \beta m.
\]
Clearly $3k_{3}\geq m;$ hence, assuming $6\beta\leq n,$
\[
12f\left(  n\right)  \beta m>\left(  n-6\beta\right)  k_{3}\left(  G\right)
\geq\left(  n-6\beta\right)  m,
\]
and the desired result follows.
\end{proof}

Applying Theorem \ref{th1} with $f\left(  n\right)  =c,$ we obtain
\begin{equation}
\sigma\left(  c\right)  \geq\frac{1}{12c+6}. \label{lob}%
\end{equation}
On the other hand, a slight modification of the graphs described in Example
\ref{examp1} gives a graph $G=G\left(  n,n^{2}/4-O\left(  1\right)  \right)
,$ such that every edge of $G$ is contained in a triangle and
\[
bk\left(  G\right)  \leq\frac{n}{6},
\]
and this, together with (\ref{lob}), implies%
\[
\lim_{c\rightarrow0}\sigma\left(  c\right)  =\frac{1}{6}.
\]
However, for large $c$ Theorem \ref{th1} is not precise enough. Prior to
obtaining a lower bound on $\gamma\left(  n,f\right)  $ that is valid in a
more general case of a function $f$, we recall the graph that Erd\H{o}s
outlined in \cite{Erd3}.

\begin{example}
Suppose $f\left(  n\right)  $ with $0<f\left(  n\right)  <n/4$ tends to
infinity with $n;$ set $l_{n}=f\left(  n\right)  ^{1/2}.$ Define a graph $G$
as follows: let $V\left(  G\right)  =\left[  n\right]  =A\cup B\cup C,$ with
$\left|  A\right|  =l_{n}^{2}$, $\left|  B\right|  =\left|  C\right|  =\left(
n-l_{n}^{2}\right)  /2.$ Join every vertex of $B$ to every vertex of $C.$
Divide $B$ and $C$ into $l_{n}$ roughly equal disjoint sets $B_{i}$ and
$C_{i}.$ Join every vertex $x_{ij}\in A$ to every vertex of $B_{i}$ and
$C_{j}.$
\end{example}

It is easily seen that $e\left(  G\right)  =n^{2}/4-f\left(  n\right)  n,$
every edge of $G$ is contained in a triangle and $bk\left(  G\right)
=o\left(  n\right)  $.

In order to obtain a precise estimate of $bk\left(  G\right)  $ we shall
describe more accurately the graph $G$. Suppose $f\left(  n\right)  $ is a
function of $n$ with $4<f\left(  n\right)  <n/4.$ Set $k=\left\lfloor \left(
2f\left(  n\right)  \right)  ^{1/2}\right\rfloor ,$ so that $k^{2}%
\leq2f\left(  n\right)  <\left(  k+1\right)  ^{2}.$ Let $n=2kt+k^{2}+s,$ where
$0\leq s<2k$. Set $V\left(  G\right)  =\left[  n\right]  $ and partition
$\left[  n\right]  $ into $2k+2$ sets $A,B_{1},...,B_{k},C_{1},...,C_{k},S$
such that
\[
\left\vert A\right\vert =k^{2},\text{ }\left\vert B_{1}\right\vert
=...=\left\vert B_{k}\right\vert =\left\vert C_{1}\right\vert =...=\left\vert
C_{k}\right\vert =t,\text{ }\left\vert S\right\vert =s.
\]

Join every vertex of $\cup_{i=1}^{k}B_{i}$ to every vertex $\cup_{i=1}%
^{k}C_{i};$ label the members of $A$ by $a_{ij}$ $\left(  i,j\in\left[
k\right]  \right)  ,$ and, for every $i,j\in\left[  k\right]  ,$ join $a_{ij}$
to all vertices of $B_{i}\cup C_{j}.$ By straightforward calculations we
obtain
\begin{align*}
e\left(  G\right)   &  =\frac{\left(  n-s-k^{2}\right)  ^{2}}{4}+k^{2}%
\frac{2\left(  n-s-k^{2}\right)  }{2k}\geq\frac{\left(  n-2k-k^{2}\right)
^{2}}{4}+k\left(  n-2k-k^{2}\right) \\
&  \geq\frac{n^{2}}{4}-\frac{k^{2}n}{2}+\frac{k^{4}-4k^{2}}{4}>\frac{n^{2}}%
{4}-f\left(  n\right)  n,
\end{align*}
and%
\[
bk\left(  G\right)  \leq\frac{n-s-k^{2}}{2k}<\frac{n}{2k}\leq\frac{n}%
{2\sqrt{2f\left(  n\right)  }}.
\]
Since, obviously, $G\in TG\left(  n,f\right)  ,$ we immediately obtain the
bound
\begin{equation}
\gamma\left(  n,f\right)  <\frac{n}{2\sqrt{2f\left(  n\right)  }}.
\label{upbs}%
\end{equation}

Our next aim is to show that, for a wide class of functions $f,$ (\ref{upbs})
is essentially tight.

\begin{theorem}
\label{th3} Let $0<c<2/5$ and $0<\varepsilon<1$ be constants, and $0<f\left(
n\right)  <n^{c}$. Then, if $n$ is sufficiently large,
\[
\gamma\left(  n,f\right)  >\left(  1-\varepsilon\right)  \frac{n}%
{2\sqrt{2f\left(  n\right)  }}.
\]

\end{theorem}

\begin{proof}
Let us start with a brief sketch of our proof\textit{.} Suppose the graph $G$
is a counterexample to our assertion. Then, from Theorem \ref{th2}, $G$ has an
induced bipartite graph $G_{1}$ of order at least $n-\alpha^{1/3}n$ and large
minimal degree. We show that each part of $G_{1}$ has cardinality close to
$n/2$ and then consider an edge from $G_{1};$ by assumption it is contained in
a triangle whose third vertex $w$ is not in $G_{1}$. We bound the degree of
$w$ from above and then bound the number of all such vertices from below.
Dropping a carefully selected number of such vertices we obtain a graph of
order $n_{1}$ and size greater than $n_{1}^{2}/4,$ such that $n_{1}$ is close
to $n.$ Then, by Corollary \ref{co1}, this graph contains a book of size
$n_{1}/6,$ completing the proof.

Now let us give the complete proof. Set $\beta=bk(G)$ and $\alpha=f\left(
n\right)  /n$\textit{.} Assume the assertion does not hold, i.e., there is
some $\varepsilon>0$ such that for every $F$ and every $N$ there is an $n>N$
with $f\left(  n\right)  >F$ and a graph $G=G\left(  n,m\right)  $ satisfying
the conditions of the theorem and with%
\begin{equation}
\beta\leq\left(  1-\varepsilon\right)  \frac{1}{2}\sqrt{\frac{n}{2\alpha}}.
\label{assum}%
\end{equation}
Then, as $\beta<n/8,$ Theorem \ref{th2} implies that $G$ has an induced
bipartite graph $G_{1}$ of order at least $n-\alpha^{1/3}n$ and
\begin{equation}
\delta\left(  G_{1}\right)  >\left(  \frac{1}{2}-4\alpha^{1/3}\right)
n=\frac{n}{2}-4\alpha^{1/3}n. \label{mindg1}%
\end{equation}
Let $V\left(  G_{1}\right)  =B\cup C$ be a bipartition of $G_{1}$ and set
$A=V\left(  G\right)  \backslash V\left(  G_{1}\right)  .$ From (\ref{mindg1}%
),
\begin{align}
\left\vert B\right\vert  &  \geq\left(  \frac{1}{2}-4\alpha^{1/3}\right)
n,\text{ }\left\vert C\right\vert \geq\left(  \frac{1}{2}-4\alpha
^{1/3}\right)  n,\label{g1ord}\\
e\left(  G_{1}\right)   &  =e\left(  B,C\right)  \geq\frac{1}{2}\left(
1-\alpha^{1/3}\right)  n\left(  \frac{1}{2}-4\alpha^{1/3}\right)  n\nonumber\\
&  =\frac{n^{2}}{4}\left(  1-\alpha^{1/3}\right)  \left(  1-8\alpha
^{1/3}\right)  >\frac{n^{2}}{4}\left(  1-9\alpha^{1/3}\right)  .\nonumber
\end{align}
Consider the set $T$ of triangles containing an edge of $G_{1}.$ Since every
edge of $G_{1}$ is contained in a triangle and $G_{1}$ is bipartite, we see
that
\begin{equation}
\left\vert T\right\vert \geq e\left(  G_{1}\right)  >\frac{n^{2}}{4}\left(
1-9\alpha^{1/3}\right)  . \label{lowT}%
\end{equation}
Let $D\subset A$ be the set of vertices of $A$ that are contained in some
triangle of $T.$ We claim that for every $w\in D,$ and $n$ sufficiently
large,
\begin{equation}
d\left(  w\right)  <\sqrt{\frac{n}{2\alpha}}. \label{updg1}%
\end{equation}
Indeed, by definition, every vertex $w\in D$ is joined to some $u\in B$ and
some $v\in C.$ Then,
\begin{align*}
\beta &  \geq\left\vert \Gamma\left(  u\right)  \cap\Gamma\left(  w\right)
\right\vert \geq\left\vert \Gamma\left(  u\right)  \cap\Gamma\left(  w\right)
\cap C\right\vert \geq d_{C}\left(  w\right)  +d_{C}\left(  u\right)
-\left\vert C\right\vert \\
&  \geq d_{C}\left(  w\right)  +\delta\left(  G_{1}\right)  -\left\vert
C\right\vert ,
\end{align*}
and, similarly,%
\[
\beta\geq\left\vert \Gamma\left(  v\right)  \cap\Gamma\left(  w\right)
\right\vert \geq\left\vert \Gamma\left(  v\right)  \cap\Gamma\left(  w\right)
\cap B\right\vert \geq d_{B}\left(  w\right)  +\delta\left(  G_{1}\right)
-\left\vert B\right\vert .
\]
Hence, summing the last two inequalities and taking into account
(\ref{mindg1}),%
\begin{align*}
2\beta &  \geq d_{B}\left(  w\right)  +d_{C}\left(  w\right)  +2\delta\left(
G_{1}\right)  -n+\left\vert A\right\vert \\
&  \geq d_{B}\left(  w\right)  +d_{C}\left(  w\right)  +\left\vert
A\right\vert -8\alpha^{1/3}n\geq d\left(  w\right)  -8\alpha^{1/3}n.
\end{align*}
To complete the proof of (\ref{updg1}), observe that from (\ref{assum}), we
have%
\[
2\beta\leq\left(  1-\varepsilon\right)  \sqrt{\frac{n}{2\alpha}}.
\]
For every $w\in D,$ let $t\left(  w\right)  $ be the number of triangles of
$T$ containing $w.$ Clearly, we have
\[
t\left(  w\right)  =\frac{1}{2}\sum_{u\in\Gamma\left(  w\right)  }\left\vert
\Gamma\left(  u\right)  \cap\Gamma\left(  w\right)  \right\vert \leq\frac
{1}{2}d\left(  w\right)  \beta\leq\frac{1}{4}d\left(  w\right)  \left(
1-\varepsilon\right)  \sqrt{\frac{n}{2\alpha}}.
\]
This, together with (\ref{updg1}), gives%
\begin{equation}
t\left(  w\right)  <\left(  1-\varepsilon\right)  \frac{n}{8\alpha}.
\label{in4}%
\end{equation}
Summing (\ref{in4}) for all $w\in D,$ in view of (\ref{lowT}), we obtain%
\[
\frac{n^{2}}{4}\left(  1-9\alpha^{1/3}\right)  <\left\vert T\right\vert
=\sum_{w\in D}t\left(  w\right)  <\left\vert D\right\vert \frac{n\left(
1-\varepsilon\right)  }{8\alpha}.
\]
Hence,%
\[
\left\vert D\right\vert >2\alpha\frac{\left(  1-9\alpha^{1/3}\right)
n}{\left(  1-\varepsilon\right)  }.
\]
Observe that, as $\alpha=f\left(  n\right)  /n<n^{c-1}$ and $c<2/5,$ we have
$\lim_{n\rightarrow\infty}\alpha^{1/3}=0.$ Then, for $n$ sufficiently large,
we see that
\[
\left\vert D\right\vert >2\left(  1+\varepsilon\right)  \alpha n.
\]
Select a set $D_{0}\subset D$ with%
\begin{equation}
\left(  2+\varepsilon\right)  \alpha n<\left\vert D_{0}\right\vert <\left(
2+2\varepsilon\right)  \alpha n. \label{bndd}%
\end{equation}
As, from (\ref{updg1}), for every vertex $w\in D_{0}$ and $n$ sufficiently
large, we have
\[
d\left(  w\right)  <\sqrt{\frac{n}{2\alpha}},
\]
then the graph $G\left[  V\backslash D_{0}\right]  $ has at least%
\[
e\left(  G\right)  -\left\vert D_{0}\right\vert \sqrt{\frac{n}{2\alpha}}%
\]
edges. We shall prove that if $n$ is large enough then
\begin{equation}
\frac{n^{2}}{4}-\alpha n^{2}-\left\vert D_{0}\right\vert \sqrt{\frac
{n}{2\alpha}}>\frac{\left(  n-\left\vert D_{0}\right\vert \right)  ^{2}}{4}.
\label{lastin}%
\end{equation}
Assume that (\ref{lastin}) does not hold. Then, from (\ref{bndd}),
\begin{align*}
\frac{n^{2}}{4}-\alpha n^{2}-\left(  2\left(  1+\varepsilon\right)  \alpha
n\right)  \sqrt{\frac{n}{2\alpha}}  &  \leq\frac{n^{2}}{4}-\alpha
n^{2}-\left\vert D_{0}\right\vert \sqrt{\frac{n}{2\alpha}}\leq\frac{\left(
n-\left\vert D_{0}\right\vert \right)  ^{2}}{4}\\
&  \leq\frac{\left(  n-\left(  2+\varepsilon\right)  \alpha n\right)  ^{2}}{4}%
\end{align*}
and thus, after some simple algebra,
\[
\frac{\varepsilon}{2}\leq\left(  2\left(  1+\varepsilon\right)  \right)
\frac{1}{\sqrt{2\alpha n}}+\frac{\left(  2+\varepsilon\right)  ^{2}\alpha^{2}%
}{4}<\frac{4}{\sqrt{2f(n)}}+4n^{2c-2},
\]
which is a contradiction if $n$ is large enough. Thus, (\ref{lastin}) holds.
Then, if $n$ is sufficiently large, Corollary \ref{co1} implies that
\[
bk\left(  G\left[  V\backslash D_{0}\right]  \right)  >\frac{n-\left\vert
D_{0}\right\vert }{6}>\sqrt{\frac{n}{2\alpha}}.
\]
This contradiction completes our proof.
\end{proof}

In \cite{Erd2}, p. 235, Erd\H{o}s asked how large $\gamma\left(
n,n^{c}\right)  $ is for $0<c<1.$ Putting $f\left(  n\right)  =n^{c-1}$ for
$1<c<7/5$ and applying Theorem \ref{th3}, together with (\ref{upbs}), we
obtain the following.

\begin{corollary}
\label{coll} If $0<c<1$ and $n$ is sufficiently large,
\[
\gamma\left(  n,n^{c}\right)  <\frac{1}{2\sqrt{2}}n^{1-c/2}.
\]
Also, if $0<c<2/5$, $\varepsilon>0$ and $n$ is sufficiently large,
\[
\ \ \ \ \ \ \ \ \ \ \ \ \ \ \ \ \ \ \ \ \ \ \ \ \ \ \ \ \ \ \ \ \ \gamma
\left(  n,n^{c}\right)  >\frac{1-\varepsilon}{2\sqrt{2}}n^{1-c/2}.
\ \ \ \ \ \ \ \ \ \ \ \ \ \ \ \ \ \ \ \ \ \ \ \ \ \ \ \ \ \ \Box
\]

\end{corollary}

\end{document}